\documentclass[a4paper,11pt]{amsart}
\usepackage{amssymb,amscd} 

\theoremstyle{plain}

\newtheorem{thm}{Theorem}[section]
\newtheorem{prop}[thm]{Proposition}
\newtheorem{lem}[thm]{Lemma}
\newtheorem{cor}[thm]{Corollary}

\theoremstyle{definition}
\newtheorem{defn}{Definition}[section]

\theoremstyle{remark}
\newtheorem{rem}{Remark}[section]
\newtheorem{exa}[rem]{Example}
\newenvironment{pf}{\medskip\noindent{Proof:}}{\hfill\qed \newline \medskip}

\title{Construction of directed strongly regular graphs using finite incidence structures}

\newcommand{\forme}[1]{}

\newcommand{\ol}{\overline}

\newcommand{\mc}{\mathcal}

\newcommand{\ts}{\textsl}
\newcommand{\vb}{\mathbf{v}}
\newcommand{\bb}{\mathbf{b}}
\newcommand{\kb}{\mathbf{k}}
\newcommand{\rb}{\mathbf{r}}
\newcommand{\lb}{\Lambda}

\begin{document} \pagenumbering{arabic} \setcounter{page}{1}
 \author[O.~Olmez]{Oktay Olmez}
\address{Department of Mathematics, Iowa State University,
Ames, Iowa, 50011, U. S. A.} \email[O. ~Olmez]{oolmez@iastate.edu}

\author[S. Y.~Song]{Sung Y. Song}
\address{Department of Mathematics, Iowa State University, Ames, Iowa, 50011, U. S.
A.} \email[S. Y.~Song]{sysong@iastate.edu}

\begin{abstract}

We use finite incident structures to construct new infinite families
of directed strongly regular graphs with parameters
\[(lq^l(q-1),\ lq^{l-1}(q-1),\ q^{l-2}(lq-l+1),\ (l-1)q^{l-2}(q-1),\
q^{l-2}(lq-l+1)),\] for integers $q$ and $l$ ($q, l\ge 2$), and
\[(lq^2(q-1),\ lq(q-1),\ lq-l+1,\ (l-1)(q-1),\ lq-l+1)\] for all prime powers $q$
and $l\in \{1, 2, \dots, q\}$. The new graphs given by these
constructions include those with parameter sets $(36, 12, 5, 2, 5)$,
$(54, 18, 7, 4, 7)$, $(72, 24, 10, 4, 10)$, $(96, 24, 7, 3, 7)$,
$(108, 36, 14, 8, 14)$ and $(108, 36, 15, 6, 15)$ listed as feasible
parameters on ``Parameters of directed strongly regular graphs" by
S. Hobart and A. E. Brouwer at ${http://homepages.cwi.nl/^\sim
aeb/math/dsrg/dsrg.html}$.  We also review these constructions and
show how our methods may be used to construct other infinite
families of directed strongly regular graphs. We then characterize
finite incidence structures from which we can produce directed
strongly regular graphs.
\end{abstract}

\maketitle

\tableofcontents


\section{Introduction and preliminaries}\label{sec:intro}

Directed strongly regular graphs were introduced by Duval \cite{Du}
in 1988 as directed versions of strongly regular graphs. There are
numerous sources for these graphs. Some of known constructions of
these graphs use combinatorial block designs \cite{FK}, coherent
algebras \cite{FK, KM}, finite geometries \cite{FK0, FK, GH, KP},
matrices \cite{Du, EH, GH}, and regular tournaments \cite{EH, Jo}.
Some infinite families of these graphs also appear as Cayley graphs
of groups \cite{DI, HS, Jo, Jo2, KM}. {\footnote{For the current
state-of-the-art on known constructions and the existence and
non-existence results of directed strongly regular graphs are found
in the Hobart and Brouwer's web \cite{BH} (also see the Handbook of
Combinatorial Designs \cite[Chap. VII.12]{BH2}.)}}

In this paper, we construct some new infinite families of directed
strongly regular graphs by using certain finite incidence
structures, such as, non-incident point-block pairs of a divisible
design and anti-flags of a partial geometry.

In Section 2, we describe a construction of a directed strongly
regular graph with parameters $(v, k, t, \lambda, \mu)$ given by
\[(lq^l(q-1),\ lq^{l-1}(q-1),\ q^{l-2}(lq-l+1),\ (l-1)q^{l-2}(q-1),\
q^{l-2}(lq-l+1)),\] for integers $l\ge 2$ and $q\ge 2$. The graph is
defined on the set of non-incident point-block pairs of group
divisible design GD$(l, q^{l-2}, q; ql)$. Among the feasible
parameters listed in \cite{BH}, our construction realizes the
feasibility of the parameter sets $(36, 12, 5, 2, 5)$ and $(96, 24,
7, 3, 7)$; and then $(72, 24, 10, 4, 10)$ and $(108, 36, 15, 6, 15)$
by applying for a construction of Duval \cite{Du}.

In section 3, we construct directed strongly regular graphs on the
sets of anti-flags of partial geometries. In particular, if we use
the partial geometry obtained from the affine plane of order $q$ by
considering all $q^2$ points and taking the $ql$ lines from exactly
$l$ pencils (parallel classes) of the plane, we obtain a directed
strongly regular graph with parameters
\[ (lq^2(q-1),\ lq(q-1),\ lq-l+1,\ (l-1)(q-1),\ lq-l+1).\]
Thus, for example, it follows the feasibility of the parameter sets
$(54, 18, 7, 4, 7)$ and $(108, 36, 14, 8, 14)$ among those listed in
\cite{BH}.

In section 4, we construct two families of directed strongly regular
graphs with parameters
\[(ql(l-1),\ q(l-1),\ q,\ 0,\ q)\] and
\[(ql(l-1),\ 2q(l-1)-1,\ ql-1,\ ql-2,\ 2q)\]
for integers $q\ge 2\mbox{ and } l\ge 3$. These graphs are not new.
J{\o}rgensen \cite{Jo2} showed that the graph with parameters
$(ql(l-1),\ q(l-1),\ q,\ 0,\ q)$ is unique for every integers $l\ge
2$ and $q\ge 2$. Godsil, Hobart and Martin constructed the second
family of graphs in \cite[Corollary 6.3]{GH}. Our construction for
the second family is essentially similar to the one in \cite{GH}.
However, both families of the graphs are constructed by using the
almost trivial incidence structure obtained from a partition of a
set of $ql$ elements into $l$ mutually disjoint $q$-element subsets.

In Section 5 and Section 6 we introduce how to construct directed
strongly regular graphs by using the anti-flags over the affine
resolvable designs and arbitrary 2-designs, respectively. In Section
7, we introduce a prolific construction method for directed strongly
regular graphs by using sets of anti-flags over certain tactical
configurations.

While there is a number of families of directed strongly regular
graphs constructed by using flags of finite geometries (cf.
\cite{BH, FK, GH, KM, KP}), there are few results for using
anti-flags or partial geometries. The results presented here adds
some constructions using sets of anti-flags of various finite
incidence structures including partial geometries, combinatorial
designs, and certain tactical configurations. More work is needed to
explore other incidence structures and address the existence and
non-existence problems of directed strongly regular graphs and their
related incidence structures.

\medskip
In the remainder of this section, we recall the definition and some
properties of directed strongly regular graphs.

A loopless directed graph $D$ with $v$ vertices is called
\textit{directed strongly regular graph} with parameters $(v, k, t,
\lambda, \mu)$ if and only if $D$ satisfies the following
conditions:
\begin{itemize}
\item[i)] Every vertex has in-degree and out-degree $k$.
\item[ii)] Every vertex $x$ has $t$ out-neighbors, all of which are
also in-neighbors of $x$.
\item[iii)] The number of directed paths of length two from a vertex
$x$ to another vertex $y$ is $\lambda$ if there is an edge from $x$
to $y$, and is $\mu$ if there is no edge from $x$ to
$y$.\end{itemize}

Another definition of a directed strongly regular graph, in terms of
its adjacency matrix, is often conveniently used. Let $D$ be a
directed graph with $v$ vertices. Let $A$ denote the adjacency
matrix of $D$, and let $I=I_v$ and $J=J_v$ denote the $v\times v$
identity matrix and all-ones matrix, respectively. Then $D$ is a
directed strongly regular graph with parameters $(v, k, t, \lambda,
\mu)$ if and only if (i) $JA=AJ=kJ$ and (ii) $A^2=tI+\lambda A+\mu
(J-I-A)$.

Duval observed that if $t=\mu$ and $A$ satisfies above equations (i)
and (ii), then so does $A\otimes J_m$ for every positive integer
$m$; and so, we have:

\begin{prop}\label{pro1} \cite{Du} If there exists a directed strongly regular
graph with parameters $(v,k,t,\lambda,\mu)$ and $t=\mu$, then for
each positive integer $m$ there exists a directed strongly regular
graph with parameters $(mv,mk,mt,m\lambda,m\mu)$.\end{prop}

\begin{prop} \label{Duval} \cite{Du} A directed strongly
regular graph with parameters $(v,k,t,\lambda,\mu)$ has three
distinct integer eigenvalues \[\theta_0=k,\quad
\theta_1=\frac12(\lambda-\mu+\delta),\quad
\theta_2=\frac12(\lambda-\mu-\delta)\] with multiplicities
\[m_0=1,\quad m_1=-\frac{k+\theta_2(v-1)}{\theta_1-\theta_2},\quad
m_2=\frac{k+\theta_1(v-1)}{\theta_1 - \theta_2},\] respectively,
where $\delta=\sqrt{(\mu-\lambda)^2+4(t-\mu)}$ is a positive
integer.\end{prop}

Throughout the paper, we will write $x\rightarrow y$ if there is an
edge from a vertex $x$ to another vertex $y$, and $x\nrightarrow y$
if there is no edge from $x$ to $y$. We will also write
$x\leftrightarrow y$ if and only if both $x\rightarrow y$ and
$y\rightarrow x$.

\section{Graphs obtained from group
divisible designs}\label{sec:CBG}

The first family of directed strongly regular graphs, which we shall
describe in this section, use non-incident point-block pairs of
group divisible designs GD$(l,q^{l-2}, q; ql)$ for integers $q\ge 2$
and $l\ge 2$.

Let $P$ be a $ql$-element set with a partition $\mc{P}$ of $P$ into
$l$ parts (`\textit{groups}') of size $q$. Let $\mc{P}=\{S_1, S_2,
\dots, S_l\}$. Let
$$\mc{B}=\{B\subset P:\ |B\cap S_i|=1 \mbox{ for all } i=1,2, \dots,
l\}.$$ Then $\mc{B}$ consists of $q^l$ subsets of size $l$ called
\textit{blocks}. The elements of $P$ will be called \textit{points}
of the incidence structure $(P, \mc{B})$ with the natural
point-block incidence relation $\in$. This structure has property
that any two points from the same group never occur together in a
block while any two points from different groups occur together in
$q^{l-2}$ blocks. It is known as a group-divisible design
GD$(l,q^{l-2}, q; ql)$.

\begin{defn} \label{D1} Let $(P, \mc{B})$ be the incidence structure
defined as above. Let $D=D(P, \mc{B})$ be the directed graph with
its vertex set
\[V(D)=\{(p, B)\in P\times \mc{B}:\ p\notin B\},\]
and directed edges given by $(p, B)\rightarrow (p', B')$ if and only
if $p\in B'$.\end{defn}

\begin{thm} \label{dsrg1} Let $D$ be the graph $D(P,\mc{B})$ defined in
Definition \ref{D1}. Then $D$ is a directed strongly regular graph
with parameters
\[\begin{array}{l} v=lq^l(q-1),\\
k= lq^{l-1}(q-1),\\ t=\mu=q^{l-2}(lq-l+1),\\
\lambda=q^{l-2} (l-1)(q-1).\\ \end{array}\]
\end{thm}

\begin{pf} It is easy to verify the values of $v$ and $k$.
 To find the value of $t$, let
$(p, B)\in V(D)$ with $B=\{b_1, b_2, \dots, b_l\}$, and let
$N^+((p,B))$ and $N^-((p,B))$ denote the set of out-neighbors and
that of in-neighbors of $(p,B)$, respectively. Then for $t$ we need
to count the elements of
\[ N^+((p,B))\cap N^-((p,B))
\ = \ \{(p^*, B^*)\in V(D):\ p\in B^* ,\ p^*\in B\}.\] Without loss
of generality, suppose that $p\in S_1$. Then condition $p\in B^*$
requires that $b_1^*$ must be $p$. The second condition $p^*\in B$
requires that $p^*=b_i$ for some $i$. With the choice of $p^*=b_1$,
there are $q^{l-1}$ blocks $B^*=\{b_1^*, b_2^*, \dots, b_l^*\}$ with
$b_1^*=p$ can be paired with $p^*=b_1$. On the other hand, with
$p^*=b_j$ for $j\neq 1$, we have $(q-1)q^{l-2}$ choices for $B^*$
with $b_1^*=p$, $b_j^*\in S_j\setminus\{b_j\}$ and $b_i^*\in S_i$
for all $i\in \{2, 3, \dots, l\}\setminus \{j\}$. Hence $t =
q^{l-1}+(l-1)(q-1)q^{l-2}$.

We now claim that the number of directed paths of length two from a
vertex $(p, B)$ to another vertex $(p', B')$ depends only on whether
there is an edge from $(p,B)$ to $(p',B')$.

Suppose that $(p, B)\rightarrow (p', B')$; that is, $p\in B'$. Then,
without loss of generality, we may assume that $p=b_1'$, and have
\[\begin{array}{ll} \lambda &=\ |N^+((p,B))\cap N^-((p',B'))|\\
& = |\{(p^*, B^*)\in V(D):\ p=b_1'\in B^*,\ p^*\in B'\}|\\
&= |\{(p^*, \{p, b_2^*,\dots, b_l^*\}): \ p^*\in
B'\setminus\{b_1'\},\ b_i^*\in S_i\setminus\{p^*\}\mbox{ for } i=2,
3, \dots, l\}|\\
& =(l-1)(q-1)q^{l-2}.\end{array} \]

For $\mu$, suppose $(p, B)\nrightarrow (p', B')$. Then
$\mu=|\{(p^*,B^*): p\in B^*,\ p^*\in B'\}|$ with $p\notin B'$.
Without loss of generality, we assume that $p\in S_1$, and have
\[\begin{array}{ll} \mu &
=\ |\{(b_1',\{p, b_2^*, \dots, b_l^*\}):\ b_i^*\in S_i \mbox{ for
}i=2, 3,\dots, l\}|\\
&\ +\ \sum\limits_{i=2}^l |\{(b_i', \{p, b_2^*,\dots, b_l^*\}):
b_i^*\in S_i\setminus\{b_i'\},\ b_j^*\in S_j\mbox{ for } j\in\{2,
3,\dots,
l\}\setminus\{i\} \}|\\
& =\ q^{l-1}+(l-1)(q-1)q^{l-2}.\end{array}\] This completes the
proof.
\end{pf}

\begin{rem} By Proposition \ref{pro1}, there are directed strongly
regular graphs with parameters
\[(mlq^l(q-1),\ mlq^{l-1}(q-1),\ mq^{l-2}(lq-l+1),\ mq^{l-2}(l-1)(q-1),\
mq^{l-2}(lq-l+1))\] for all positive integers $m$. These graphs can
be constructed directly by replacing all edges by multiple edges
with multiplicity $m$ in the above construction.
\end{rem}
\begin{rem} In the case of $l=2$ ($m=1$) in the above construction, we
obtain the graphs with parameters
$$(2q^2(q-1), \ 2q(q-1),\ 2q-1,\ q-1,\ 2q-1).$$ These graphs are
constructed on the sets of non-incident vertex-edge pairs of
complete bipartite graphs $K_{q,q}$. If we use complete bipartite
multigraph with multiplicity $m$ for each edge as in the above
remark, we obtain directed strongly regular graphs with parameters,
$$(2mq^2(q-1), \ 2mq(q-1),\ m(2q-1),\ m(q-1),\ m(2q-1)).$$ With
several combinations of small $m$ and $q$, we obtain new directed
strongly regular graphs with parameter sets $(36, 12, 5, 2, 5)$,
$(96, 24, 7, 3, 7)$, $(72, 24, 10, 4, 10)$ and $(108, 36, 15, 6,
15)$. So the feasibility of these parameter sets which are listed in
\cite{BH} has been realized.
\end{rem}

By Proposition \ref{Duval}, we can easily compute the eigenvalues
for $D(P,\mc{B})$.
\begin{lem}\label{cor1} If a directed strongly
regular graph with parameters $(v,k,t,\lambda,\mu)$ satisfies
$t=\mu=\lambda +q^{l-1}$ for positive integers $q$ and $l\ge 2$,
then its eigenvalues are $\theta_0=k,\ \theta_1=0,\
\theta_2=-q^{l-1}$ with multiplicities $m_0=1,\
m_1=v-1-\frac{k}{q^{l-1}},\ m_2=\frac{k}{q^{l-1}},$ respectively.
\end{lem}

\begin{exa} (Parameter sets for small DSRGs (orders up to 110)
constructed.)
\[\begin{array}{ll}
(v, k, t, \lambda, \mu) & l, q, m\\
\hline
(8, 4, 3, 1, 3) & l=q=2\\
(16, 8, 6, 2, 6) &  l=q=2;\ m=2\\
(24, 12, 9, 3, 9) &  l=q=2;\ m=3\\
(24, 12, 8, 4, 8) & l=3, q=2\\
(32, 16, 12, 4, 12) &  l=q=2;\ m=4\\
(36, 12, 5, 2, 5) & l=2, q=3 \\ 
(40, 20, 15, 5, 15)&  l=q=2;\ m=5\\
(48, 24, 16, 8, 16) &  l=3, q=2;\ m=2\\
(48, 24, 18, 6, 18) & l=q=2;\ m=6\\
(56, 28, 21, 7, 21) & l=q=2;\ m=7\\
(64, 32, 20, 4, 20) & l=4, q=2\\
(64, 32, 24, 8, 24) & l=q=2;\ m=8\\
(72, 24, 10, 4, 10) &  l=2, q=3;\ m=2\\
(72, 36, 24, 12, 24) & l=3, q=2; \ m=2\\
(72, 36, 27, 9, 27) &  l=q=2; \ m=9\\
(80, 40, 30, 10, 30) &  l=q=2;\ m=10\\
(88, 44, 33, 11, 33) & l=q=2;\ m=11\\
(96, 24, 7, 3, 7)& l=2, q=4\\
(96, 48, 32, 16, 32) &  l=3, q=2;\ m=4\\
(96, 48, 36, 12, 36) &  l=q=2;\ m=12\\
(104, 52, 39, 13, 39) & l=q=2;\ m=13\\
(108, 36, 15, 6, 15) &  l=2, q=3;\ m=3\\ \end{array}\]
\end{exa}

\section{Graphs obtained from partial geometries}\label{sec:AP}

In this section, we use partial geometries to construct a family of
directed strongly regular graphs. The concept of a partial geometry
was introduced by R. C. Bose in connection with his study of large
cliques of more general strongly regular graphs in \cite{Bo}.

A partial geometry $\ts{pg}(\kappa, \rho, \tau)$ is a set of points
$P$, a set of lines $\mc{L}$, and an incidence relation between $P$
and $\mc{L}$ with the following properties:
\begin{enumerate}
\item Every line is incident with $\kappa$ points ($\kappa\ge 2$),
and every point is incident with $\rho$ lines ($\rho\ge 2$).
\item Any two points are incident with at most one line.
\item If a point $p$ and a line $L$ are not incident, then
there exists exactly $\tau$ ($\tau\ge 1$) lines that are incident
with $p$ and incident with $L$.\end{enumerate} Here we use
parameters $(\kappa, \rho, \tau)$ instead of more traditional
notations, $(K, R, T)$ or $(1+t, 1+s, \alpha)$ used in \cite{Bo, BC}
or \cite{Th}. In what follows, we often identify a line $L$ as the
set of $\kappa$ points that are incident with $L$; so, when $p$ and
$L$ are incident, we write ``$p\in L$," as well as ``$p$ is on $L$,"
and ``$L$ passes through $p$".

\begin{defn}\label{dpg} Let $D=D(\ts{pg}(\kappa, \rho, \tau))$ be the directed
graph with its vertex set
\[V(D)=\{(p, L)\in P\times \mc{L}:\ p\notin L\},\]
and directed edges given by $(p, L)\rightarrow (p', L')$ if and only
if $p\in L'$.
\end{defn}

\begin{thm}\label{pg} Let $D$ be the graph $D(\ts{pg}(\kappa, \rho,\tau))$
defined as above. Then $D$ is a directed strongly regular graph with
parameters
\[\begin{array}{l} v=\frac{\kappa\rho(\kappa-1)(\rho-1)}{\tau}
\left(1+\frac{(\kappa-1)(\rho-1)}{\tau}\right),\\
k=\frac{\kappa\rho(\kappa-1)(\rho-1)}{\tau},\\
t=\mu= \kappa\rho-\tau,\\
\lambda= (\kappa-1)(\rho-1).\\ \end{array}\]
\end{thm}

\begin{pf} The value $v$ is clear as it counts the number of anti-flags
of $\ts{pg}(\kappa,\rho,\tau)$ which has
$\kappa\left(1+\frac{(\kappa-1)(\rho-1)}{\tau}\right)$ points and
$\rho\left(1+\frac{(\kappa-1)(\rho-1)}{\tau}\right)$ lines. Also it
is clear that \[k=|\{(p',L'): p\in
L'\}|=\rho(v-\kappa)=\kappa\rho(\kappa-1)(\rho-1)/\tau.\]

For $t$, given a vertex $(p, L)\in V(D)$, we count the cardinality
of
\[
N^+((p,L))\ \cap\ N^-((p,L))=\{(p', L')\in V(D):\ p\in L' ,\ p'\in
L\}.
\] Among the $\rho$ lines passing through $p$, $\rho-\tau$ lines
are parallel to $L$. If $L'$ is parallel to $L$, then all $\kappa$
points on $L$ can make legitimate non-incident point-line pairs
$(p', L')$ with given $L'$. In the case when $L'$ is not parallel to
$L$, all points on $L$ except for the common incident point of $L$
and $L'$, can form desired pairs. Hence we have
\[t=(\rho-\tau)\kappa + \tau (\kappa-1)= \kappa\rho-\tau.\]

For $\lambda$, suppose $(p, L)\rightarrow (p',L')$. It is clear that
 \[\lambda=|\{(p^*,L^*): p\in L^*,\ p^*\in L'\}|=(\rho-1)(\kappa-1),\]
 because each of $\rho-1$ lines passing through $p$ (excluding $L'$)
 can be paired with any of $\kappa-1$ points on $L'\setminus\{p\}$.

For given $(p, L)$ and $(p',L')$ with $p\notin L'$,
 \[\mu=|\{(p^*,L^*): p\in L^*,\ p^*\in L'\}|=(\rho-\tau)\kappa
 +\tau(\kappa-1)=\kappa\rho-\tau,\]
since among the $\rho$ lines passing through $p$, the ones that are
parallel to $L'$ can form desired pairs with any of $\kappa$ points
on $L'$, while each of the remaining $\tau$ lines can be paired with
$\kappa-1$ points on $L'$.
\end{pf}

Partial geometries are ubiquitous. For example, partial geometries
with $\tau =1$ are generalized quadrangles, those with $\tau=\kappa
-1$ are transversal designs, and those with $\tau=\rho-1$ are known
as nets. In order to have some concrete examples of new directed
strongly regular graphs of small order, we consider a special class
of partial geometries that are obtained from finite affine planes.

Let $AP(q)$ denote the affine plane of order $q$. Let $\ol{AP}^l(q)$
denote the partial geometry obtained from $AP(q)$ by considering all
$q^2$ points and taking the lines of $l$ parallel classes of the
plane. Then $\ol{AP}^l(q)$ inherits the following properties from
$AP(q)$: (i) every point is incident with $l$ lines, and every line
is incident with $q$ points, (ii) any two points are incident with
at most one line, (iii) if $p$ and $L$ are non-incident point-line
pair, there are exactly $l-1$ lines containing $p$ which meet $L$.
That is, $\ol{AP}^l(q)=\ts{pg}(q,l,l-1)$.

\begin{cor} \label{AP} Let $D=D(\ol{AP}^l(q))$ be the directed graph
$D(\ts{pg}(q, l, l-1))$ defined as in Definition \ref{dpg}. Then $D$
is a directed strongly regular graph with parameters
\[(v,k, t, \lambda, \mu)=(lq^2(q-1),\ lq(q-1),\ lq-l+1,\ (l-1)(q-1),\
lq-l+1).\] \end{cor}

\begin{pf} It immediately follows from Theorem \ref{pg}. \end{pf}

In particular, if $l=q$, $\ol{AP}^q(q)=\ts{pg}(q,q,q-1)$ is a
transversal design TD$(q,q)=(P, \mathcal{G}, \mathcal{L})$ of order
$q$, block size $q$, and index 1 in the following sense.
\begin{enumerate}
\item  $P$ is the set of $q^2$ points of $\ol{AP}^q(q)$;
\item  $\mc{G}$ is the partition of $P$ into $q$ classes (groups) such that
each class consists of $q$ points that were collinear in $AP(q)$ but
not in $\ol{AP}^q(q)$;
\item  $\mc{L}$ is the set of $q^2$ lines (blocks);
\item every unordered pair of points in $P$ is contained in either
exactly one group or in exactly one block, but not both.
\end{enumerate}

\begin{cor} Let $D$ be the graph $D(\ol{AP}^q(q))$ defined as above.
Then $D$ is a directed strongly regular graph with parameters
\[(v,k, t, \lambda, \mu)=(q^3(q-1),\ q^2(q-1),\ q^2-q+1,\ (q-1)^2,\
q^2-q+1).\] The eigenvalues of this graph are $q^2(q-1),\ 0,\ -q$
with multiplicities $1,\ q^4-q^3-q^2+q-1,\ q(q-1),$ respectively.
\end{cor}

\begin{exa} This method produces the DSRGs with the following
parameter sets (with order up to $110$ and $m=1$ only):
\[\begin{array}{ll}
(v, k, t, \lambda, \mu) & l, q, m\\
\hline (8, 4, 3, 1, 3) & l=q=2\\
(12, 6, 4, 2, 4) & l=3, q=2\\
(16, 8, 5, 3, 5) & l=4, q=2\\
(20, 10, 7, 4, 7) & l=5, q=2\\
(24, 12, 7, 5, 7) & l=6, q=2\\
(28, 14, 8, 6, 8) & l=7, q=2\\
(32, 16, 9, 7, 9) & l=8, q=2\\
(36, 12, 5, 2, 5) & l=2, q=3 \\
(54, 18, 7, 4, 7) & l=q=3\\
(72, 24, 9, 6, 9) & l=4, q=3\\
(96, 24, 7, 3, 7) & l=2, q=4\\
\end{array}\]
Among these, the parameter sets $(54, 18, 7, 4, 7)$; and so $(108,
36, 14, 8, 14),$ $(162, 54, 21, 12, 21) \cdots $ confirm the
feasibility of the putative parameter sets listed in \cite{BH}.
These graphs share the same automorphism group which is isomorphic
to $(((Z_3\times Z_3)\rtimes Z_3)\rtimes Z_2)\rtimes Z_2$.
\end{exa}

\begin{rem} When $l=2$, we obtain the directed strongly regular graph
with parameters $(4q^2, 4q, 2q-1, q-1, 2q-1)$ in both Section 2 and
Section 3. The two graphs are shown to be isomorphic although the
construction methods are different. We demonstrate the isomorphism
through an example; namely, for the case when $q=3$, the directed
strongly regular graph with parameters $(36, 12, 5, 2, 5)$. The one,
denoted by $D_1$, is coming from non-incident vertex-edge pairs of
the complete bipartite graph $K_{3,3}$ with partite sets $\{1, 2,
3\}$ and $\{4, 5, 6\}$. The other, denoted by $D_2$, comes as
$D_2=D(\overline{AP})^2(3)$, with point set $P=\{1, 2, \dots, 9\}$
and line set $\mc{L}=\{123, 456, 789, 147, 258, 369\}$ of
$\overline{AP}^2(3)$.

The following map between the vertices (anti-flags of the underlying
incident structures,) establishes the isomorphism between the graphs
$D_1$ and $D_2$, where the adjacency of vertices in $D_1$ is defined
by
$$(h, ij)\rightarrow (h', i'j')\mbox{ iff } h\in i'j'$$ while that
of $D_2$ is defined by
$$(hij,l)\rightarrow (h'i'j',l')\mbox{ iff } l'\in \{hij\}.$$

\[\begin{array}{|lcl|lcl|lcl|}\hline
1,24 &\leftrightarrow & 123,4 & 2,14 &\leftrightarrow & 456, 1 & 3,14 &\leftrightarrow & 789, 1\\
1,25 &\leftrightarrow & 123,5 & 2,15 &\leftrightarrow & 456, 2 & 3,15 &\leftrightarrow & 789, 2\\
1,26 &\leftrightarrow & 123,6 & 2,16 &\leftrightarrow & 456, 3 & 3,16 &\leftrightarrow & 789, 3\\
1,34 &\leftrightarrow & 123,7 & 2,34 &\leftrightarrow & 456, 7 & 3,24 &\leftrightarrow & 789, 4\\
1,35 &\leftrightarrow & 123,8 & 2,35 &\leftrightarrow & 456, 8 & 3,25 &\leftrightarrow & 789, 5\\
1,36 &\leftrightarrow & 123,9 & 2,36 &\leftrightarrow & 456, 9 & 3,26 &\leftrightarrow & 789, 6\\
4,15 &\leftrightarrow & 147,2 & 5,14 &\leftrightarrow & 258, 1 & 6,14 &\leftrightarrow & 369, 1\\
4,16 &\leftrightarrow & 147,3 & 5,16 &\leftrightarrow & 258, 3 & 6,15 &\leftrightarrow & 369, 2\\
4,25 &\leftrightarrow & 147,5 & 5,24 &\leftrightarrow & 258, 4 & 6,24 &\leftrightarrow & 369, 4\\
4,26 &\leftrightarrow & 147,6 & 5,26 &\leftrightarrow & 258, 6 & 6,25 &\leftrightarrow & 369, 5\\
4,35 &\leftrightarrow & 147,8 & 5,34 &\leftrightarrow & 258, 7 & 6,34 &\leftrightarrow & 369, 7\\
4,36 &\leftrightarrow & 147,9 & 5,36 &\leftrightarrow & 258, 9 & 6,35 &\leftrightarrow & 369, 8\\
\hline
\end{array}\]

\end{rem}

\section{Graphs obtained from partitioned sets}\label{sec:DSGN}

In this section we construct directed strongly regular graphs for
two parameter sets,
\[(ql(l-1),\ q(l-1),\ q,\ 0,\ q)\] and \[(ql(l-1),\ 2q(l-1)-1,\
ql-1,\ ql-2,\ 2q)\] for all positive integers $q\ge 1$ and $l\ge 3$.
The incidence structure which will be used here may be viewed as a
degenerate case of those used in earlier sections. The graphs
produced by this construction are not new. The latter family of
graphs have been constructed by Godsil, Hobart and Martin in
\cite[Corollary 6.3]{GH}. Nevertheless, we introduce our
construction to illustrate a variation of constructions which may be
applied to produce new families of directed strongly regular graphs.

Let $P$ be a set of $ql$ elements (`points'), and let $S_1, S_2,
\dots, S_l$ be $l$ mutually disjoint $q$-element subsets of $P$
(`blocks'). We denote the family of blocks by $\mc{S}=\{S_1, S_2,
\dots, S_l\}$. We will say that point $x\in P$ and block $S\in
\mc{S}$ is a non-incident point-block pair if and only if $x\notin
S$.

For the directed strongly regular graph with the first parameter
set, let $D=D(P,\mc{S})$ be the directed graph with its vertex set
\[V(D)=\{(x, S)\in P\times \mc{S}:\ x\notin
S \},\] and directed edges defined by $(x,S)\rightarrow (x',S')$ if
and only if $x\in S'$.

\begin{thm} Let $P$, $\mc{S}$ and $D(P,\mc{S})$ be as the above.
Then $D(P,\mc{S})$ is a directed strongly regular graph with
parameters
\[(ql(l-1),\ q(l-1),\ q,\ 0,\ q).\]
\end{thm}
\begin{pf}
Easy counting arguments give the values for parameters.
\end{pf}

\begin{rem} J{\o}rgensen \cite{Jo2} showed that the directed strongly regular
graph with parameters $(ql(l-1),\ q(l-1),\ q,\ 0,\ q)$ is unique for
all positive integers $l$ and $q$.
\end{rem}

For the directed strongly regular graph with the second parameter
set, let $G=G(P,\mc{S})$ be the directed graph with its vertex set
\[V(G)=\{(x, S)\in P\times \mc{S}:\ x\notin S \},\] where $P$ and
$\mc{S}$ are as the above, and let edges on $V(G)$ be defined by:
$(x,S)\rightarrow (x',S')$ if and only if
\[ \left \{\begin{array}{l} (1)\ x\in S' ; \mbox{ or} \\
(2)\ S=S'\mbox{ and } x\neq x'.\\ \end{array}\right . \]

\begin{thm} Let $G$ be the graph $G(P,\mc{S})$ defined as above.
Then $G$ is a directed
strongly regular graph with parameters \[(ql(l-1),\ 2q(l-1)-1,\
ql-1,\ ql-2,\ 2q).\]\end{thm}
\begin{pf} Clearly $v=ql(l-1)$ as $|V(G)|=|V|\cdot (|\mc{S}|-1)$.

Given a vertex $(x,S)$, let \[\begin{array}{l}
N_1((x, S)):=\{(x', S'): x\in S' \mbox{ and } x'\in S\}\\
N_2((x, S)):=\{(x', S'): S=S' \mbox{ and }
x\neq x'\}\\
N_3((x, S)):=\{(x', S'): x\in S'\mbox{ and } x'\notin S\}\\
\end{array}\]
Then by simple counting, we have
\[\begin{array}{l} |N_1((x, S))|=q\\ |N_2((x, S))|= q(l-2)+q-1\\
|N_3((x, S)|=q(l-2). \\ \end{array}\] Hence, $t=|N_1((x, S))|+
|N_2((x, S))|= ql-1$, and $k=t+|N_3((x,S))|=2ql-2q-1$.

Finally, we count the number of vertices $(x^*, S^*)$ that belongs
to $N^+((x,S))\cap N^-((x',S'))$ for $\lambda$ and $\mu$.

For $\lambda$, given an edge $(x,S)\rightarrow (x',S')$, we shall
need to consider the following two cases separately: (Case 1) when
$x\in S'$, and (Case 2) when $S=S'$ and $x\neq x'$.\\ Case (1).
Suppose $x\in S'$. Then unless $S^*=S$, $S^*$ must contain $x$ which
forces $S^*=S'$. With $S^*=S'$, all points that are not belong to
$S'\cup \{x'\}$ can be chosen to be $x^*$; and so, there are
$q(l-1)-1$ choices for $x^*$. With $S^*=S$, the $q-1$ points of the
set $S'\setminus\{x\}$ are possible for $x^*$. Together there are
$ql-2$
vertices $(x^*, S^*)$ belonging to $N^+((x,S))\cap N^-((x',S'))$.\\
Case (2). Let $S=S'$. Then with the choice of $S^*=S=S'$, $x^*$ can
be chosen from $P\setminus (S\cup \{x, x'\})$; and so, there are
$q(l-1)-2$ choices for $x^*$. Also with the choice of $S^*$ being
the block containing $x$, any point in $S$ can be chosen as $x^*$;
and so, there are $q$ possible choices for $x^*$. Together we have
$ql-2$ choices for $(x^*, S^*)$ as well. Hence we have
$\lambda=ql-2$.

For $\mu$, suppose $S\neq S'$ and $x\notin S'$. Then it is clear
that
\[\mu=|\{(x^*, S^*): x\in S^*, x^*\in S'\}\cup \{(x^*, S^*): S^*=S,
x\in S'\}|=q+q.\] This completes the proof.
\end{pf}
\section{Graphs obtained from affine resolvable designs}
Let $\vb, \kb$ and $\lb$ be positive integers such that
$\vb>\kb>\lb$. Let $\mc{D}=(P, \mc{B})$ be a $2-(\vb,\kb,\lb)$
design such that $|P|=\vb$, each block consists of $\kb$ points, and
each pair of distinct points occurs in exactly $\lb$ blocks. Every
point belongs to $\rb=\frac{\lb(\vb-1)}{\kb-1}$ blocks, and the
block set $\mc{B}$ consists of $\bb=\frac{\lb\vb(\vb -1)}{\kb(\kb
-1)}$ blocks. A parallel class is a set of pairwise disjoint blocks
that is a partition of $P$. Obviously a parallel class can exist
only if $\vb$ is divisible by $\kb$ or equivalently, $\bb$ is
divisible by $\rb$.

A $2-(\vb,\kb,\lb)$ design is resolvable if the block set $\mc{B}$
can be partitioned into $\rb$ disjoint parallel classes. A
resolvable design $\mc{D}=(P, \mc{B})$ is called affine (or affine
resolvable) if there exists a positive integer $m$ such that any two
non-parallel blocks intersect in exactly $m$ point(s). All
parameters of an affine design may be expressed in terms of $m$ and
another parameter $s$ where $s=({\vb}/{\kb})=({\bb}/{\rb})$, the
number of blocks in a parallel class, as follows (cf. \cite{Bose}):
\[\vb=m s^2,\ \bb=\frac{s(m s^2-1)}{s-1},\ \rb=\frac{m s^2 -1}{s
-1},\ \kb=m s,\ \lb=\frac{m s -1}{s-1}.\]

Let $(P, \mc{B})$ be an affine resolvable $2-(\vb, \bb, \kb, \rb,
\lb)$ design with $m={\kb^2}/{\vb}$ and $s={\vb}/{\kb}$. We denote
this design by $AP_{m}(s)$. Let $\ol{AP}^l_{m}(s)$ denote the
incidence structure obtained from $AP_{m}(s)$ by considering all
$\vb$ points in $P$ and taking the blocks of any $l$ parallel
classes of the design $(P, \mc{B})$ for $1<l< \rb$. Let
$\bar{\mc{B}}$ denote the set of blocks in the $l$ parallel classes.
Then $\ol{AP}^l_{m}(s)=(P, \bar{\mc{B}})$ inherits the following
properties from $AP_{m}(s)$: (i) every point is incident with $l$
blocks, and every block is incident with $\kb$ points, (ii) any two
points are incident with at most $\lb$ blocks, (iii) if $p$ and $B$
are non-incident point-block pair, there are exactly $l-1$ blocks
containing $p$ which intersect with $B$ in $m$ point(s).

\begin{thm} \label{dsrg-ARD} For $m\ge 1, s\ge 2$ and $l\ge 2$,
let $\ol{AP}^l_{m}(s)$ be the tactical configuration defined as
above. Let $D=D(AP^l_{m}(s))$ be the directed graph with its vertex
set
\[V(D)=\{(p, B)\in P\times \bar{\mc{B}}:\ p\notin B\},\]
and directed edges defined by $$(p, B)\rightarrow (p', B') \mbox{ if
and only if } p\in B'.$$ Then $D$ is a directed strongly regular
graph with parameters
\[v=m ls^2(s-1),\ k=m ls(s-1),\ t=\mu=m(ls-l+1),\ \lambda= m(l-1)(s-1).\]
\end{thm}
\begin{pf} By straightforward counting, we find
\[v=ls(\vb -\kb),\ k=l(\vb -\kb),\ t=(l-1)(\kb -m)+\kb,\
\lambda= (l-1)(\kb -m),\ \mu=(l-1)(\kb -m)+\kb.\] If we write these
parameters in terms of $l, m$ and $s$, we get the desired
expression.
\end{pf}

\begin{rem}
This is a generalization of the construction of the directed
strongly regular graphs obtained from the degenerated affine plane
described in Corollary \ref{AP}. In fact, the graphs obtained for
$m=1$ in Theorem \ref{dsrg-ARD} are the ones obtained in Corollary
\ref{AP}. From the graphs with $m=1$, we can obtain graphs with
$m\ge 2$ by the `$m$-multiple construction method' introduced by
Duval \cite{Du}.
\end{rem}

\begin{rem}
In regard to the source of affine resolvable designs, we would like
to recall the concept of orthogonal arrays. An orthogonal array
$A[n, q, \mu]$ of strength 2, $q$ constraints and index $\mu$ is a
$q\times \mu n^2$ matrix with entries from an $n$-set with the
property that any two of its rows contain each of the $n^2$ ordered
pairs of symbols exactly $\mu$ times. The parameters of an
orthogonal array satisfy the inequality $q\le \frac{\mu
n^2-1}{n-1}$. An orthogonal array is called complete if if satisfy
the equality $q= \frac{\mu n^2-1}{n-1}$. The existence of an affine
resolvable $2-(\vb,\bb, \kb, \rb, \Lambda)$ is equivalent to the
existence of a complete orthogonal array with $n=\frac{\vb}{\kb}$,
$q=\rb$, $\mu=\frac{\kb^2}{\vb}$ \cite{Sh}. Many concrete examples
may be obtained from the study of orthogonal arrays. \end{rem}

Here are the parameters of the graphs (with small orders) that are
constructed by our method described above. -- This table may be
deleted.

\[\begin{array} {|rrrrr|l||rrrrr|l|}\hline
v &  k & t& \lambda & \mu & (\vb, \kb, \lb; s, m, l) & v &  k & t&
\lambda & \mu & (\vb, \kb, \lb; s, m, l)\\
\hline \hline 12 & 6& 4& 2 & 4 &  (4, 2,6; 2, 1, 3)
&132& 66& 36& 30& 36& (12,6,22;2,3,11) \\
\hline 8 & 4& 3& 1& 3& (4,2,6; 2, 1, 2)
&120& 60& 33& 27& 33& (12,6,22;2,3,10)\\
\hline 56& 28& 16& 12& 16& (8,4,14;2, 2, 7)&
108& 54& 30& 24& 30& (12,6,22;2,3,9)\\
\hline 48 & 24& 14& 10& 14& (8,4,14;2, 2, 6)
&96& 48& 27& 21& 27& (12,6,22;2,3,8)\\
\hline 40 & 20& 12& 8& 12& (8,4,14; 2, 2,5)
&84& 42& 24& 18& 24& (12,6,22;2,3,7)\\
\hline 32 & 16& 10& 6& 10& (8,4,14;2,2,4)
&72& 36& 21& 15& 21& (12,6,22;2,3,6)\\
\hline 24 & 12& 8& 4& 8& (8,4,14;2,2,3)
&60& 30& 18& 12& 18& (12,6,22;2,3,5)\\
\hline 16 & 8& 6& 2& 6& (8,4,14;2, 2, 2)
&48& 24& 15& 9& 15& (12,6,22;2,3,4)\\
\hline 72 & 24& 9& 6& 9& (9,3,12;3,1,4)
&36& 18& 12& 6& 12& (12,6,22;2,3,3)\\
\hline *54& 18& 7& 4& 7& (9,3,12;3,1,3)
&24& 12& 9& 3& 9& (12,6,22;2,3,2)\\
\hline *36& 12& 5& 3& 5& (9,3,12;3,1,2) & & & & & &\\
 \hline \end{array}\]

\bigskip

\section{Graphs obtained from $2-(\vb, \kb, \lb)$ designs}

The balanced incomplete block designs have been used to produce
directed strongly regular graphs (cf. \cite{BH}). In this section,
we will see that we can also construct many directed strongly
regular graphs on the sets of anti-flags of every $2-(\vb, \bb, \kb,
\rb, \lb)$ design (with $\bb+\lb>2\rb$) in the following two
different ways.

\begin{thm} \label{dsrg-2-D1}
Let $\mc{D}=(P, \mc{B})$ be a $2-(\vb, \bb, \kb, \rb, \lb)$ design
with $\bb+\lb>2\rb$.
 Let $D=D(P,\mc{B})$ be the directed graph with its vertex set
\[V(D)=\{(p, B)\in P\times {\mc{B}}:\ p\notin B\},\]
and directed edges defined by $(p, B)\rightarrow (p', B')$ if and
only if $p'\in B$.

Then $D$ is a directed strongly regular graph with parameters
\[v=\vb (\bb - \rb)=\frac{\lb \vb (\vb -1)(\vb - \kb)}
{\kb (\kb -1)},\quad k=\kb (\bb -\rb)=\frac{\lb (\vb -1)(\vb -
\kb)}{\kb (\kb -1)},\]
\[t=\mu=\kb (\rb -\lb)=\frac{\lb \kb (\vb - \kb)}{\kb -1},\quad
\lambda= (\kb -1)(\rb -\lb)=\lb(\vb - \kb).\]
\end{thm}

\begin{pf} It is easy to verify the values of the parameters.
\end{pf}

\begin{rem} This is another way to construct the graphs that may be
constructed by the method of [T8] in Brouwer and Hobart's home page.
The construction [T8] uses the complementary designs of 2-designs
instead of 2-designs themselves.
\end{rem}

The second method of obtaining directed strongly regular graphs from
the anti-flags of 2-designs uses a different rule for the adjacency.
\begin{thm} \label{dsrg-2-D2}
Let $\mc{D}=(P, \mc{B})$ be a $2-(\vb, \bb, \kb, \rb, \lb)$ design
with $\bb+\lb>2\rb$.
 Let $D=D(P,\mc{B})$ be the directed graph with its vertex set
\[V(D)=\{(p, B)\in P\times {\mc{B}}:\ p\notin B\},\]
and directed edges given by $(p, B)\rightarrow (p', B')$ if and only
if either $p'\in B$ or $p=p'$ and $B\neq B'$. Then $D$ is a directed
strongly regular graph with parameters
\[\begin{array}{l} v=\vb (\bb - \rb),\\
k=\kb (\bb -\rb)+(\bb -\rb -1),\\ t=\kb (\rb -\lb)+ (\bb -\rb -1),\\
\lambda= \kb (\rb -\lb) + (\bb -\rb -2), \\ \mu=(\kb +1) (\rb -\lb).\\
\end{array}\]
\end{thm}

\begin{pf} The values of the parameters can be verified by the similar
counting argument.
\end{pf}


\begin{thebibliography}{99}


\bibitem{Bose} R. C. Bose. A note on the resolvability f incomplete
block designs, \textit{Sankhy$\tilde{a}$} \textbf{6} (1942),
105--110.

\bibitem{Bo} R. C. Bose. Strongly regular graphs, partial
geometries, and partially balanced designs, \textit{Pacific J.
Math.} \textbf{13} (1963) 389--419.

\bibitem{BC} A. E. Brouwer, A. M. Cohen and A. Neumaier.
\textit{Distance-Regular Graphs}, Springer-Verlag, Berlin, 1989.


\bibitem{BH} A. E. Brouwer and S. Hobart. Parameters of directed strongly regular graphs.
http://homepages.cwi.nl /\~ aeb/math/dsrg/dsrg.html

\bibitem{BH2} A. E. Brouwer and S. Hobart. Directed strongly regular graph,
in: C. Colburn and J. Dinitz (Eds.), \textit{Handbook of
Combinatorial Designs}, CRC Inc., Boca Raton, 868--875.

\bibitem{CV} C. J. Colburn and S. A. Vanstone. Doubly resolvable
twofold triple systems, \textit{Congr. Numer.} 34(1982) 219-233.

\bibitem{Du} A. Duval. A Directed Graph Version of Strongly Regular Graphs.
{\em Journal of Combinatorial Theory}, Series A \textbf{47} (1988),
71--100.

\bibitem{DI} A. Duval and D. Iourinski. Semidirect product constructions
of directed strongly regular graphs, {\em Journal of Combinatorial
Theory}, (A) \textbf{104} (2003) 157-167.

\bibitem{EH} C. Eschenbach, F. Hall, R. Hemasinha, S. Kirkland, Z. Li,
B. Shader, J. Stuart, J. Weaver. Properties of Tournaments among
Well-Matched Players. {\em The American Mathematical Monthly},
\textbf{107} (2000), 881--892.

\bibitem{FK0} F. Fiedler, M. Klin and M. Muzychuk. Small
vertex-transitive directed strongly regular graphs, {\em Discrete
Mathematics} \textbf{255} (2002) 87-115.

\bibitem{FK} F. Fiedler, M. Klin, Ch. Pech. Directed strongly
regular graphs as elements of coherent algebras, in: K. Denecke,
H.-J. Vogel (Eds.), General Algebra and Discrete Mathematics, Shaker
Verlag, Aachen, 1999, pp. 69--87.

\bibitem{Gr} M. Greig. Log tables and block designs. \textit{Bull.
Inst. Combin. Appl.} 48 (2006), 66–72.

\bibitem{GH} C. D. Godsil, S. A. Hobart and W. J. Martin.
Representations of directed strongly regular graphs. {\em European
J. Combin.} \textbf{28} (2007), no. 7, 1980--1993.

\bibitem{HS} S. Hobart, and T. Shaw.  A note on a family of directed strongly
regular graphs. {\em European Journal of Combinatorics} 20 (1999),
819--820.

\bibitem{Jo1} L. J{\o}rgensen. Search for directed strongly regular graphs.
Report R-99-2016   Department of Mathematical Sciences, Aalborg
University, 1999.

\bibitem{Jo} L. J{\o}rgensen. Directed strongly regular graphs with
$\mu$ = $\lambda$. {\em Discrete Mathematics}, \textbf{231} (2001),
no. 1--3, 289--293.

\bibitem{Jo2} L. J{\o}rgensen. Non-existence of directed Strongly Regular Graphs.
{\em Discrete Mathematics}, \textbf{264} (2003), 111--126.




\bibitem{KM} M. Klin, A. Munemasa, M. Muzychuk and P.-H. Zieschang.
Directed strongly regular graphs obtained from coherent algebras,
{\em Linear Algebra Appl.} \textbf{377} (2004), 83--109.

\bibitem{KP} M. Klin, C. Pech and P.-H. Zieschang. Flag algebras
of block designs: I. Initial notions, Steiner 2-designs, and generalized
quadrangles, preprint MATH-AL-10-1998, Technische Universitat Dresden.

\bibitem{SS} K. See and S. Y. Song. Association schemes of small
order, {\em Journal of Statistical Planning and Inference}
\textbf{73} (1998) 225--271.

\bibitem{Sh} S. S. Shrikhande. Affine resolvable balanced
incomplete block designs: a survey, \textit{Aequationes Math.}
\textbf{14} (1976), 251--269.

\bibitem{Th} J. A. Thas. Partial geometries,
in: C. Colburn and J. Dinitz (Eds.), Handbook of Combinatorial
Designs, CRC Inc., Boca Raton, 557--561.

\end{thebibliography}
\end{document}